\documentclass{amsart}

\usepackage{amssymb, array, verbatim, amscd, xypic}
\usepackage{amsmath, amsfonts,amsthm}

\theoremstyle{plain}
\newtheorem{theorem}{Theorem}[section]
\newtheorem{corollary}[theorem]{Corollary}
\newtheorem{lemma}[theorem]{Lemma}
\newtheorem{proposition}[theorem]{Proposition}

\theoremstyle{definition}

\newtheorem{example}[theorem]{Example}

\newtheorem{remark}[theorem]{Remark}

\newcommand{\CF}{{\mathcal F}}

\newcommand{\boldP}{{\mathbf P}}

\newcommand{\Ima}{\mathrm{Im}}

\newcommand{\Ker}{\mathrm{Ker}}

\newcommand{\Tor}{\mathrm{Tor}}

\newcommand{\frakF}{\mathfrak{F}}

\newcommand{\Hom}{\operatorname{Hom}}

\newcommand{\Ext}{\operatorname{Ext}}

\begin{document}

\title[Modules $M$ such that $\Ext_R^1(M,-)$ commutes with direct limits]{Modules $M$ such that $\Ext_R^1(M,-)$ commutes with direct limits}
\author{Simion Breaz}
\thanks{Research supported by the
CNCS-UEFISCDI grant PN-II-RU-TE-2011-3-0065}
\address{Babe\c s-Bolyai University, Faculty of Mathematics
and Computer Science, Str. Mihail Kog\u alniceanu 1, 400084
Cluj-Napoca, Romania}
\email{bodo@math.ubbcluj.ro}

\date{\today}


\thanks{ }

\begin{abstract}
We will use Watts's theorem together with Lenzing's
characterization of finitely presented modules via commuting
properties of the induced tensor functor in order to study
commuting properties of covariant Ext-functors.

 \end{abstract}


\maketitle
\date{}

\section{Introduction}

It is well known that commuting properties of some canonical
functors (as Hom or tensor functors induced by a right module)
provide important information (about that module) or some
important tools in the study of some subcategories for the module
category. For instance, H. Lenzing proved in \cite[Satz 3]{Le}
that a right $R$-module $M$ is finitely presented if and only if
the functor $\Hom_R(M,-)$ preserves direct limits (i.e. filtered
colimits) or the tensor product $-\otimes_R M$ commutes with
direct products, \cite[Satz 2]{Le}. These theorems had a great
influence in modern algebra: the first result is used to define
finitely presented objects in various categories, e.g. \cite{AR},
while the second is an important ingredient in Chase's
characterization of right coherent rings \cite[Theorem
2.1]{Chase}.

The property that the covariant Hom-functor commutes with direct
sums provides us with the notion of \textsl{small} (compact)
module, \cite[p.54]{Ba}. This notion is useful in many topics in
module theory as generalizations of Morita equivalences
\cite{CoF}, almost free modules \cite{Tr95} or the internal
structure of the ring, \cite{ZT}. It is well known that every
finitely generated module (these are the modules such that the
induced covariant Hom-functor commutes with direct unions,
\cite[24.10]{Wis}) is small. Moreover, for some important classes
of rings (as right noetherian or right perfect by \cite{Re},
\cite{CT}) these two conditions are equivalent, but there are also
important types of rings for which there are non-finitely
generated small modules (e.g. non-artinian regular simple rings),
\cite{EGT}. It is proved in \cite{Zem-cluj} that every small
module is finitely presented if and only if the ring is
noetherian. A similar problem can be proposed for the covariant
functor $\Ext^1_R$: \textit{Identify classes of rings $R$ such
that a functor $\Ext^1_R(M,-)$ commutes with direct unions (or
direct limits) if and only if it commutes with direct sums.}
Corollary \ref{cor-hereditary} provides an answer to this problem.
Moreover, in Example \ref{example-com} it is shown that these two
conditions are not equivalent in general. We mention here that $R$
is right coherent exactly if for every right $R$-module $M$ the
functor $\Ext_R^1(M,-)$ commutes with direct unions if and only if
it commutes with direct limits, \cite[Corollary 7]{Br12}.

Concerning commuting properties of the derived functors $\Ext^*$
and $\Tor^*$, it was proved by Brown in \cite{Brown} that
Lenzing's results can be extended to produce more finiteness
conditions. In order to state this results, let us recall from
\cite[p.~103]{GT06} that a right $R$ module $M$ \textsl{is
$FP_n$}, for a fixed integer $n\geq 0$, if $M$ has a projective
resolution which is finitely generated in every dimension $\leq
n$.

\begin{theorem}\label{th-brown} \cite[Theorem 1]{Brown}, \cite[Theorem A]{Strebel}
The following are equivalent for a right $R$ module $M$ and a
non-negative integer $n$:
\begin{enumerate}
\item $M$ is $FP_{n+1}$;

\item $\Ext_R^k(M,-)$ commutes with direct limits for all $0\leq
k\leq n$,

\item $\Tor^R_k(M,-)$ commutes with direct products for all $0\leq
k\leq n$.
\end{enumerate}
\end{theorem}

The proof of this theorem is based on Lenzing's characterizations
of finitely presented modules via commuting properties of
covariant Hom-functors with direct limits, respectively commuting
properties of tensor functors with direct products, and these are
applied to obtain independently (1)$\Leftrightarrow$(2) and
(1)$\Leftrightarrow$(3). For the case of $\Ext$-functors, this
theorem was refined and completed by Strebel in \cite{Strebel}. We
can ask if there is a cyclic proof. A connection is suggested in
\cite[Corollary 1]{Strebel}, where it is proved that a right
$R$-module $M$ of projective dimension at most $n$ has the
property that $\Ext_R^n(M,-)$ commutes with direct limits if and
only if $M$ has a projective resolution which is finitely
generated in dimension $n$.

In the present paper we give such a proof in Theorem
\ref{theorem-fp-omega}, where modules of projective dimension at
most 1 such that the induced covariant $\Ext^1$-functor commutes
with direct limits are characterized via Lenzing's
characterization of finitely presented modules by commuting
properties of the tensor product. This theorem is used to obtain
some results from Strebel's paper, \cite{Strebel}. Similar
techniques were used by Krause in \cite{Krause98} in order to
characterize general coherent functors.


In the end of this introduction let us observe that there are
commuting properties which are non-trivial for the Hom-functors
but they are trivial for $\Ext$-functors. For instance, it is well
known that the covariant Hom-functors commute with inverse limits.
However, as in \cite[Example 3.1.8]{GT06}, we can use
\cite[Theorem 2]{Be11} to write every module as an inverse limit
of injective modules. Therefore

\begin{proposition}
Let $M$ be a module. The functor $\Ext_R^1(M,-)$ commutes with
inverse limits if and only if $M$ is projective.
\end{proposition}

Dualizing the definition of small modules we obtain the notion of
slender module. The structure of these modules can be very
complicated (see \cite{EM02} or \cite{GT06}). Transferring this
approach to the contravariant Ext-functor, we say that the
cotravariant functor $\Ext_R^n(-,M)$ \textsl{inverts products} if
for every countable family $\CF=(M_i)_{i<\omega}$ of right
$R$-modules the canonical homomorphism
$\bigoplus_{i<\omega}\Ext_R^n(M_i,M)\to \Ext_R^n(\prod_{i<\omega}
M_i,M)$ is an isomorphism. However, in the case when $M$ has the
injective dimension at most $n$ this property is very restrictive.

\begin{proposition}\label{ext-slender}
Let $M$ be a right $R$-module of injective dimension at most
$n\geq 1$. The functor $\Ext^n_R(-,M)$ inverts products with
countable many factors if and only if $M$ is of injective
dimension at most $n-1$.
\end{proposition}

\begin{proof}
By the dimension shifting formula, it is enough to assume $n=1$.

In this hypothesis, we start with a module $N$ and with the
canonical monomorphism $0\to N^{(\omega)}\to N^\omega$. Applying
the contravariant Ext-functor and \cite[Proposition 2.4]{ABS} we
obtain the natural homomorphism
$$\Ext^1_R(N,M)^{(\omega)}\cong\Ext^1_R(N^\omega,M)\to
\Ext^1_R(N^{(\omega)},M)\cong \Ext^1_R(N,M)^\omega,$$ which
coincides with the canonical homomorphism
$$\Ext^1_R(N,M)^{(\omega)}\to \Ext^1_R(N,M)^\omega,$$ and it is an
epimorphism. This is possible only if $\Ext^1_R(N,M)=0$.
\end{proof}

We do not know what happens in Proposition \ref{ext-slender} in
case we do not assume any bound on the injective dimension. It is
also an open question when the contravariant $\Ext^1$-functor
preserves products. In \cite{GP} the authors provide an answer for
abelian groups, and we refer to \cite{Br11} for the case of
contravariant Hom-functors. Other commuting properties of these
functors, for the case of abelian groups, are studied in
\cite{ABS12} and \cite{S11}.

\section{When $\Ext$ commutes with direct limits}

Let $R$ be a unital ring and $M$ a right $R$-module. If
$\frakF=(M_{i}, \upsilon_{ij})_{i,j\in I}$ is a direct system of
modules and $\upsilon_{i}:M_i\to \underrightarrow{\lim}M_i$ are
the canonical homomorphisms, then for every non-negative integer
$n$ there is a canonical homomorphism
$$\Phi^{n,M}_\frakF:\underrightarrow{\lim}\Ext^n_R(M,M_i)\to \Ext^n_R(M,\underrightarrow{\lim}
M_i),$$ the natural homomorphism induced by the family
$\Ext^n_R(M,\upsilon_{ij})$, $i,j\in I$.

For the case of direct sums, if we consider $\CF=(M_i)_{i\in I}$ a
family of modules and we denote by $u_i:M_i\to \bigoplus_{i\in I}
M_i$ the canonical homomorphisms then we will obtain a natural
homomorphism
$$\Phi^{n,M}_\CF:\bigoplus_{i\in I}\Ext^n_R(M,M_i)\to \Ext^n_R(M,\bigoplus_{i\in I}
M_i),$$  induced by the family $\Ext^n_R(M,u_i)$, $i\in I$.

We will use these homomorphisms for the cases $n\in\{0,1\}$, and
we will denote $\Phi^{0,M}_\frakF:=\Psi^{M}_\frakF$,
$\Phi^{1,M}_\frakF=\Phi^{M}_\frakF$, respectively
$\Phi^{0,M}_\CF:=\Psi^{M}_\CF$, $\Phi^{1,M}_\CF=\Phi^{M}_\CF$.

We say that $\Ext^n_R(M,-)$ \textsl{commutes with direct limits}
(\textsl{direct sums}) if the homomorphisms $\Phi^{n,M}_\frakF$
(respectively $\Phi^{n,M}_\CF$) are isomorphisms for all directed
systems $\frakF$ (respectively all families $\CF$).


Lenzing's theorem says us that $M$ is finitely presented if and
only if $\Psi^M_\frakF$ are isomorphisms for all $\frakF$.
Moreover, using Theorem \ref{th-brown} for $n=1$, we observe that
$M$ is an $FP_2$ module if and only if $\Psi^M_\frakF$ and
$\Phi^M_\frakF$ are isomorphisms for all families $\frakF$. We
will start with a slight improvement of this result, replacing the
condition ``$\Psi^M_\frakF$ is an isomorphism'' (i.e. $M$ is
finitely presented) by the condition ``M is finitely generated''.
Therefore in the case of finitely generated modules the hypothesis
of \cite[Lemma 3.1.6]{GT06} is sharp.

Following the terminology used in \cite{Br12}, we will call a
module $M$ an \textsl{fp-$\Omega^1$-module}
(\textsl{fg-$\Omega^1$-module}, respectively
\textsl{small-$\Omega^1$-module}) if there is a projective
resolution
$$(\boldP):\ \dots\to P_2\to P_1\overset{\alpha_1}\to P_0\to M\to 0$$ such
that the first syzygy $\Omega^1(\boldP)=\Ima(\alpha_1)$ is
finitely presented (finitely generated, respectively small).

\begin{lemma}\label{lemma-fp-omega}
If $M$ is fp-$\Omega^1$ (fg-$\Omega^1$, respectively
small-$\Omega^1$) right $R$-modules then there is a projective
module $L$ such that $M\oplus L$ is a direct sum of an
$FP_2$-module (finitely presented module, respectively finitely
generated module) and a projective module.
\end{lemma}

\begin{proof}
Let $0\to K\overset{\alpha}\to P\overset{\beta}\to M\to 0$ be an
exact sequence such that $P$ is projective and $K$ is finitely
presented (finitely generated, respectively small). If $C$ is a
projective module such that $P\oplus C\cong R^{(I)}$ is free, then
we consider the induced exact sequence $0\to K\overset{\alpha}\to
P\oplus C\overset{\beta\oplus 1_C}\to M\oplus C\to 0$. Since $K$
is small (recall that every finitely generated module is small) as
a right $R$-module there is a finite subset $J$ of $I$ such that
$\Ima(\alpha)\subseteq R^{(J)}$. Then $M\oplus C\cong
R^{(J)}/\Ima(\alpha)\oplus R^{(I\setminus J)}=H\oplus L$, where
$H$ is an $FP_2$-module (finitely presented, respectively finitely
generated) and $L$ is projective.
\end{proof}

The closure under direct summands of the class of
fg-$\Omega^1$-modules can be characterized by the commuting
property of the covariant $\Ext^1$-functor with direct unions
(i.e. direct system of
monomorphisms). 
We state this result for sake
of completeness.

\begin{theorem} \cite[Theorem 5]{Br12}
The following are equivalent for a right $R$-module $M$:
\begin{enumerate}
\item $M$ is a direct summand of an fg-$\Omega^1$-module;


\item The functor $\Ext_R^1(M,-)$ commutes with direct systems of
monomorphisms.
\end{enumerate}
\end{theorem}

For the cases of fp-$\Omega^1$-modules and
small-$\Omega^1$-modules we are able to prove a similar result
only in the case of finitely generated modules. In fact, it is
easy to observe, using \cite[Lemma 3.1.6]{GT06} that the covariant
$\Ext^1$-functor induced by an fp-$\Omega^1$-modules commutes with
direct limits. The converse is true for finitely generated
modules, \cite[Lemma 1]{Br12}. A similar result is valid for
small-$\Omega^1$-modules.

\begin{theorem}\label{t1}
Let $M$ be a right $R$-module $M$.
\begin{enumerate}
\item If $M$ is an fp-$\Omega^1$-module then $\Ext^1_R(M,-)$
commutes with direct limits.


\item If $M$ is a small-$\Omega^1$-module then $\Ext^1_R(M,-)$
commutes with direct sums.
\end{enumerate}

The converses of these statements are true if $M$ is finitely
generated.
\end{theorem}

\begin{proof}
We will prove (2). The proof for (1) follows the same steps. Using
Lemma \ref{lemma-fp-omega}, we can suppose that $M$ is finitely
generated.

If $0\to K\to P\to M\to 0$ is an exact sequence with $P$ a
projective module then for every family
 $\CF=(M_{i})_{i\in I}$ of right $R$-modules, we have the
following useful commutative diagram
$$(\sharp)\ \ \  \begin{CD}
0@>>> \oplus_{i\in I}\Hom_R(M,M_i) @>>> \oplus_{i\in I}\Hom_R(P,M_i) @>>> \  \\
@.
@VV{\Psi^M_\CF}V @VV{\Psi^P_\CF}V @. \\
0@>>> \Hom_R(M,\oplus_{i\in I} M_i) @>>>
\oplus_{i\in I} \Hom_R(P,M_i) @>>> \ \\
\ @>>> \oplus_{i\in I}\Hom_R(K,M_i) @>>>
\oplus_{i\in I}\Ext^1_R(M,M_i) @>>> 0 \\
@. @VV{\Psi^K_\CF}V  @VV{\Phi^M_\CF}V @.
\\ \ @>>>\Hom_R(K,\oplus_{i\in I} M_i)@>>> \Ext^1_R(M,\oplus_{i\in I} M_i) @>>> 0 \end{CD}$$
whose rows are exact.

Since $M$ and $P$ are finitely generated, the arrows
$\Psi^M_\frakF$ and $\Psi^P_\frakF$ are isomorphisms by
\cite[24.10]{Wis}. Using Five Lemma we observe that
$\Psi^K_\frakF$ is an isomorphism if and only if $\Phi^M_\frakF$
is an isomorphism.
\end{proof}

We are ready to state the main result of this paper. It states
that fp-$\Omega^1$-modules of projective dimension at most 1 can
be characterized by commuting properties of the induced covariant
$\Ext^1$-functor.

\begin{theorem}\label{theorem-fp-omega}
The following are equivalent for a right $R$-module $M$ of
projective dimension 1:
\begin{enumerate}

\item  $M$ is an fp-$\Omega^1$-module.

\item There is a projective module $L$ such that $M\oplus L$ is a
direct sum of an $FP_2$-module ($FP_1$-module) and a projective
module.

\item $\Ext_R^1(M,-)$ commutes with direct limits.

\item $\Ext_R^1(M,-)$ commutes with direct sums of copies of $R$.
\end{enumerate}
\end{theorem}

\begin{proof}
We only need to prove $(4)\Rightarrow(1)$.

Let $M$ be a module of projective dimension at most 1 such that
$\Ext^1_R(M,-)$ commutes with direct sums of copies of $R$.
Therefore $\Ext^1_R(M,-)$ is a right exact functor which commutes
with direct sums of copies of $R$. By Watts's theorem
\cite[Theorem 1]{Watts} (and its proof) we conclude that the
functor $\Ext^1_R(M,-)$ is naturally equivalent to the functor
$-\otimes_R \Ext^1_R(M,R)$. It follows that the tensor product
functor $-\otimes_R \Ext^1_R(M,R)$ preserves the products, hence
that $\Ext^1_R(M,R)$ is a finitely presented left $R$-module.

Let $n$ be a positive integer such that the left $R$-module
$\Ext^1_R(M,R)$ is generated by $n$ elements. Using \cite[Lemma
6.9]{Tr-w} we conclude that there is an exact sequence of right
modules $0\to R^n\to C\to M\to 0$ such that $\Ext^1_R(C,R)=0$.
Starting with this exact sequence we obtain for every cardinal
$\kappa$ a commutative diagram
$$\begin{CD}
\Hom(R^n,R)^{(\kappa)}@>>> \Ext^1_R(M,R)^{(\kappa)}@>>>
\Ext^1_R(C,R)^{(\kappa)} @>>> 0 \\
@VVV @VVV @VVV @. \\ \Hom(R^n,R^{(\kappa)})@>>>
\Ext^1_R(M,R^{(\kappa)})@>>> \Ext^1_R(C,R^{(\kappa)}) @>>> 0
\end{CD},$$
where the vertical arrows are the canonical homomorphisms induced
by the universal property of direct sums. Moreover the first
vertical is an isomorphism since every finitely generated module
is small, while the second is also an isomorphism by our
hypothesis. Therefore the third vertical arrow is also an
isomorphism. It follows that $\Ext^1_R(C,R^{(\kappa)})=0$ for all
cardinals $\kappa$. Since $C$ is of projective dimension at most
1, it is not hard to see that $C$ is projective.
\end{proof}




\begin{corollary}\label{cor-hereditary}
Let $R$ be an right hereditary ring and $M$ a right $R$-module.
The following are equivalent:
\begin{enumerate}
\item $\Ext^1_R(M,-)$ commutes with direct limits;

\item $\Ext^1_R(M,-)$ commutes with direct sums;

\item $\Ext^1_R(M,-)$ commutes with direct sums of copies of $R$;

\item $M=N\oplus P$, where $N$ is finitely presented and $P$ is
projective.
\end{enumerate}
\end{corollary}

\begin{proof}
Since every right hereditary ring is right coherent, every
finitely presented module is an $FP_2$-module. Therefore only
$(3)\Rightarrow(4)$ requires a proof.

Suppose that $\Ext_R^1(M,-)$ commutes with direct sums of copies
of $R$. Using Theorem \ref{theorem-fp-omega}, we observe that
there is a projective module $L$ such that $M\oplus L= F\oplus U$
with $F$ finitely presented and $U$ projective. The conclusion
follows using the same techniques as in \cite{Al61}. If
$\pi_U:F\oplus U\to U$ is the canonical projection, then
$\pi_U(M)$ is projective. Then $F\cap M=\Ker(\pi_{U|M})$ is a
direct summand of $M$. Therefore $N=F\cap M$ is a direct summand
of $F\oplus U=M\oplus L$. Using this, it is not hard to see that
$N=F\cap M$ is a direct summand of $F$. Hence $N$ is finitely
presented and $M=N\oplus P$, where $P\cong \pi_U(M)$ is
projective.
\end{proof}

\begin{example}\label{example-com}

The equivalence $(3)\Leftrightarrow(4)$ from Theorem
\ref{theorem-fp-omega} is not valid for general finitely generated
modules.

To see this, it is enough to consider a ring $R$ which has a
non-finitely generated small right ideal $I$. For instance, it is
proved in \cite{Tr} that every direct product of infinitely many
rings has such an ideal. Since $R/I$ is a finitely generated right
$R$-module, we can apply Theorem \ref{t1} to see that
$\Ext^1_R(R/I,-)$ commutes with direct sums, but it does not
commute with direct limits.
\end{example}



\begin{remark}
Corollary \ref{cor-hereditary} was proved for abelian groups in
\cite{BS11} without the condition that the isomorphisms
$\bigoplus_{i\in I}\Ext^1_R(M,M_i)\cong \Ext^1_R(M,\bigoplus_{i\in
I} M_i)$ are the natural ones.
\end{remark}

\begin{corollary} \cite[Corollary 1]{Strebel}
Let $M$ be right $R$-module of finite projective dimension $\leq
n$. Then $M$ admits a projective resolution which is finitely
generated in dimension $n$ if and only if $\Ext_R^n(M,-)$
preserves direct sums of copies of $R$. If one of the equivalent
conditions is fulfilled, then the functors $\Ext_R^n(M,-)$ and
$-\otimes_R \Ext_R^n(M,R)$ are naturally isomorphic.
\end{corollary}

\begin{proof}
Let $(\mathbf{P})$ be a projective resolution of $M$. Then the
$(n-1)$-th syzygy $\Omega^{n-1}(\mathbf{P})$ is of projective
dimension at most 1. Using the dimension shifting formula we
observe that $\Ext_R^n(M,-)$ is naturally isomorphic to
$\Ext_R^1(\Omega^{n-1}(\mathbf{P}),-)$, and the conclusion follows
from Theorem \ref{theorem-fp-omega} and Watts's theorem.
\end{proof}

\begin{lemma}\label{lemma-quotient}
Let $L$ be a right $R$-module of projective dimension at most $1$
and $K$ a finitely presented submodule of $L$. The module $M=L/K$
has the property that $\Ext_R^1(M,-)$ preserves direct sums of
copies of $R$ if and only if $M$ is an fp-$\Omega^1$-module.
\end{lemma}

\begin{proof}
Let $I$ be a set. It induces a commutative diagram
$$\begin{CD} \Hom_R(K,R)^{(I)}@>>>\Ext_R^1(M,R)^{(I)}@>>>
\Ext_R^1(L,R)^{(I)}@>>> \\ @VV{\Psi_K}V @VV{\Phi_M}V @VV{\Phi_L}V \\
\Hom_R(K,R^{(I)})@>>>\Ext_R^1(M,R^{(I)})@>>>
\Ext_R^1(L,R^{(I)})@>>>\\ \ @>>> \Ext_R^1(K,R)^{(I)} @>>>
\Ext_R^2(M,R)^{(I)}@>>> 0 \\ @. @VV{\Phi_K}V @VVV @. \\ \ @>>>
\Ext_R^1(K,R^{(I)})@>>> \Ext_R^2(M,R^{(I)})@>>> 0
\end{CD}$$
where all arrows represent the canonical homomorphisms. Since $K$
is finitely presented, the arrows $\Psi_K$ and $\Phi_K$ are
isomorphisms (for $\Phi_K$ we use Theorem \ref{t1}). Moreover,
$\Phi_M$ is an isomorphism from our hypothesis and it is obvious
that the last vertical arrow is also a monomorphism by \cite[Lemma
2.2]{Strebel}. Therefore $\Ext_R^1(L,-)$ commutes with direct sums
of copies of $R$. By Theorem \ref{theorem-fp-omega}, it follows
that there is a projective resolution $0\to R^n\to P\to L\to 0$.

Using all these data we construct a pullback diagram $$\begin{CD}
@. @. 0 @. 0 @. \\ @. @. @VVV @VVV @.\\
0@>>> R^n@>>> U@>>> K @>>> 0 \\ @. @| @VVV @VVV @. \\
0@>>> R^n@>>> P@>>> L @>>> 0 \\ @. @. @VVV @VVV @.\\ @. @. M@= M
@.\\ @. @. @VVV @VVV @.\\ @. @. 0@. 0 @.
\end{CD}.$$
Since $U$ is finitely presented, the proof is complete.
\end{proof}

We also obtain the principal ingredient for the proof of
\cite[Theorem B]{Strebel}.

\begin{corollary} \label{strebel2.6}\cite[Lemma 2.6]{Strebel}
Let $M$ be a right $R$-module such that $\Ext_R^1(M,-)$ commutes
with direct sums of copies of $R$. Suppose that there is a
projective resolution $$ (\boldP)\ \ P_2 \to
P_1\overset{\alpha_1}\to P_0\overset{\alpha_0}\to M\to 0$$ such
that $\Omega^2(\boldP)$ is finitely generated. Then there is a
projective resolution
$$(\mathbf{P'})\ \ P_2\to P'_1\to P'_0\to
M\to 0$$ with $P'_1$ finitely generated such that
$\Omega^2(\boldP)=\Omega^2(\mathbf{P'})$.
\end{corollary}

\begin{proof}
We can suppose that there is a set $I$ such that $P_1\cong
R^{(I)}$. If $J$ is a finite subset of $I$ such that
$\Omega^2(\boldP)\subseteq R^{(J)}$ and
$K=R^{(J)}/\Omega^2(\boldP)$ then we have an exact sequence $$0\to
K\oplus R^{(I\setminus J)}\overset{\alpha}\to P_0\to M\to 0,$$
where $\alpha$ is induced by $\alpha_1$. If
$M'=P_0/\alpha(R^{(I\setminus J)})$ then we have an exact sequence
$0\to K\to M'\to M\to 0$. Observe that $M'$ is of projective
dimension at most 1. By Lemma \ref{lemma-quotient}, the module $M$
is an fp-$\Omega^1$-module. Therefore there is an exact sequence
$0\to U\to Q\to M\to 0$ with $Q$ projective and $U$ finitely
presented, and we can construct a commutative diagram with exact
rows
$$\begin{CD}
@. 0 @. 0 @. 0 @. \\ @. @VVV @VVV @VVV @. \\
0@>>> \Omega^2(\boldP) @>>> V @>>> U @>>> 0 \\ @. @VVV @VVV @VVV
@. \\ 0@>>> P_1 @>>> W@>>> Q @>>> 0\\
@. @VVV @VVV @VVV @.\\
0@>>> \Omega^1(\boldP)@>>> P_0@>>> M @>>> 0 \\
@. @VVV @VVV @VVV @. \\ @. 0 @. 0 @. 0 @.
\end{CD}$$
Since the middle vertical row is splitting, $V$ is a finitely
generated projective module. Moreover we can construct an exact
sequence $0\to \Omega^2(\boldP)\to V\to Q\to M\to 0$, and the
proof is complete.
\end{proof}



\noindent{\textbf{Acknowldgements:}} I thank Jan Trlifaj for
drawing my attention on Strebel's work \cite{Strebel}. I also
thank to the Referee, whose comments improved this paper.


 \end{document}